\pdfoutput=1

\documentclass[12pt]{amsart} 
\usepackage{graphicx,color}
\usepackage[margin=2cm]{geometry}
\usepackage{amsthm,amsmath,amssymb,enumerate,comment}
\usepackage{mathtools}

\usepackage{cleveref}
\usepackage{autonum}
\usepackage[all]{xy}
\usepackage{here}
\usepackage{tikz}
\usetikzlibrary{calc}
\usetikzlibrary {arrows.meta}

\newcommand{\F}{\mathbb{F}}

\newcommand{\Q}{\mathbb{Q}}
\newcommand{\R}{\mathbb{R}}

\newcommand{\Z}{\mathbb{Z}}

\newcommand{\confi}[2]{C_{#1}({#2})}
\newcommand{\uconfi}[2]{B_{#1}({#2})}
\newcommand{\torus}[1]{T^{#1}}
\newcommand{\sphere}[1]{S^{#1}}
\newcommand{\rp}[1]{\mathbb{R}P^{#1}}
\newcommand{\sym}[1]{\Sigma_{#1}}

\newcommand{\chomology}[3]{H_{#1}({#2};{#3})}

\newcommand{\ccohomology}[3]{H^{#1}({#2};{#3})}
\newcommand{\rccohomology}[3]{\tilde{H}^{#1}({#2};{#3})}
\newcommand{\lcohomology}[3]{\mathcal{H}^{#1}({#2};{#3})}

\newcommand{\abhtpygrp}[2]{\pi_{#1}({#2})}
\newcommand{\cellchain}[2]{C_{#1}({#2})}
\newcommand{\cellcochain}[3]{C^{#1}({#2};{#3})}

\newcommand{\map}[3]{{#1}\colon{#2}\rightarrow {#3}}
\newcommand{\namelessmap}[2]{{#1}\rightarrow {#2}}

\newcommand{\iset}[2]{\{{#1}|{#2}\}}

\newcommand{\Hom}[3]{\text{Hom}_{#1}({#2},{#3})}

\renewcommand{\leq}{\leqslant}
\renewcommand{\geq}{\geqslant}

\renewcommand{\hat}{\widehat}
\renewcommand{\tilde}{\widetilde}

\newtheoremstyle{mystyle}
    {3pt}
    {3pt}
    {\itshape}
    {}
    {\bfseries}
    {}
    {10pt}
    {\thmname{#1} \thmnumber{#2}\thmnote{(#3)}}

\theoremstyle{mystyle}

\newtheorem{theorem}{Theorem}[section]

\newtheorem{lemma}{Lemma}[section]

\numberwithin{equation}{subsection}

\makeatletter
\renewenvironment{proof}[1][\proofname]{\par
  \pushQED{\qed}%
  \normalfont \topsep6\p@\@plus6\p@\relax
  \trivlist
  \item\relax
  {#1\@addpunct{.}}\hspace\labelsep\ignorespaces
}{%
  \popQED\endtrivlist\@endpefalse
}
\makeatother

\renewcommand{\proofname}{\textbf{Proof}}

\crefname{section}{Chapter}{Chps.}
\crefname{corollary}{Corollary}{Cors.}
\crefname{equation}{eq.}{eqs.}
\crefname{align}{eq.}{eqs.}
\crefname{figure}{Fig.}{figs.}
\crefname{theorem}{Theorem}{Theorems}
\crefname{lemma}{Lemma}{Lems.}

\title{Mod 2 Representation of the Symmetric Group of Order 2 over Cohomology Groups of 2-Configuration Space of Torus}
\author{Tomoki TOKUDA}
\email{tokuda.tomoki.379@s.kyushu-u.ac.jp}
\address{Joint Graduate School of Mathematics for Innovation, Kyushu University, Fukuoka, 819-0395, Japan}

\begin{document}

\maketitle
\markboth{TOMOKI TOKUDA}{MOD 2 REPRESENTATION OF $\sym{2}$}

\begin{abstract}
   In this article, we compute the mod 2 representarion of the symmetric group of order 2 over the singular cohomology groups of orderd 2-configuration space $C_{2}(T^{d})$ of the $d$-torus $T^{d}$ for $d\geq 1$. As applications of the computation, we determine the Stiefel-Whitney height of $C_{2}(T^{d})$ for any $d$, and determine $H^{*}(\mathbb{R}P^{\infty};\mathbb{F}_2)$-module structure of the cohomology groups of the unordered 2-configuration space of the $d$-torus for $d=2,3$ using the Serre spectral sequence.
\end{abstract}

\section{Introduction and Main Theorems}
The $n$-configuration space $\confi{n}{X}$ of a space $X$ is the subspace of the $n$-fold product space $X^{n}$ defined by 
\begin{align}
    \confi{n}{X}=\iset{(x_1,\dots ,x_n)\in X^{n}}{x_i\neq x_j \text{ if } i\neq j}.
\end{align}
The symmetric group $\sym{n}$ acts on $\confi{n}{X}$ freely by permuting coordinates. The unordered $n$-configuration space is the orbit space $\uconfi{n}{X}=\confi{n}{X}/\sym{n}$. In general, given an action of a group $G$ on a space $X$, the singular chain complex of $X$ inherits the $G$-action. Then $\chomology{i}{X}{R}$ and $\ccohomology{i}{X}{R}$ become a module over the group ring $R[G]$ where $R$ is a ring with unity. 

\vspace{3mm}
The homology groups of the unordered configuration spaces $\uconfi{n}{M}$ of a manifold $M$ are studied by Bödigheimer, Cohen and Taylor in \cite{BODIGHEIMER1989111}. And this result is generalised by Chen and Zhang in \cite{chen2022mod}. By their works, the Betti numbers of $\chomology{i}{\uconfi{n}{M}}{k}$ are almost determined where $k$ is a prime field. 

The cohomology algebra of unordered configuration spaces over closed manifolds was computed by Félix and Tanré in \cite{10.1112/S0024610705006794} with coefficients $\Q$ and $\F_p$. For the $\uconfi{n}{\torus{2}}$, Pagaria computed its cohomology ring structure in \cite{10.2140/agt.2020.20.2995}. However, in all these papers, the characteristics of the coefficient fields are required to be zero or greater than the number of coordinates of configurations. 

\vspace{3mm}
The main result of this paper is the description of the $\F_2[\sym{2}]$-module structure of the cohomology groups $\ccohomology{i}{\confi{2}{\torus{d}}}{\F_2}$ for $d\geq 2$. The point of this study is that the characteristic of the field $\F_2$ divides the order of $\sym{2}$. In general, any representation of a finite group $G$ over a finite dimensional $k$-vector space $V$ is completely reducible unless the characteristic of $k$ divides the order of $G$. However, without the assumption, the representation become more complicated in general.

\vspace{3mm}
Here is the main result of this paper.
\begin{theorem}\label{confi_rep}
    The representation of $\sym{2}$ over $\ccohomology{i}{\confi{2}{\torus{d}}}{\F_2}$ is decomposed as following:
    \begin{align}
        \ccohomology{i}{\confi{2}{\torus{d}}}{\F_2}\cong
        \begin{cases}
            \F_2^{\oplus \binom{d}{k}}\oplus \F_2[\sym{2}]^{\oplus \frac{1}{2}\left(\sum_{j=0}^{i}\binom{d}{j}\binom{d}{i-j}-\binom{d}{k}\right)}  & \text{  if $i=2k<d$,}\\
            \F_2[\sym{2}]^{\oplus \frac{1}{2}\sum_{j=0}^{i}\binom{d}{j}\binom{d}{i-j}}  & \text{  if $i=2k+1<d$,}\\
            \F_2^{\oplus\binom{d}{k}+\binom{d}{i-d} }\oplus \F_2[\sym{2}]^{\oplus \frac{1}{2}\left(\sum_{j=0}^{i}\binom{d}{j}\binom{d}{i-j}-\binom{d}{k}\right)-\binom{d}{i-d}}  & \text{  if $d\leq i=2k<2d$ ,}\\
            \F_2^{\oplus\binom{d}{i-d} }\oplus \F_2[\sym{2}]^{\oplus \frac{1}{2}\left(\sum_{j=0}^{i}\binom{d}{j}\binom{d}{i-j}-\binom{d}{k}\right)-\binom{d}{i-d}}  & \text{  if $d\leq i=2k+1<2d$ ,}\\
            0  & \text{  if $i\geq 2d$.}\\
        \end{cases}
    \end{align}
    
\end{theorem}

\vspace{5mm}
As a one of the additional results, we can determine the $\F_2[\sym{2}]$-module structure of $\ccohomology{i}{\uconfi{2}{\torus{d}}}{\F_2}$ for $d=2,3$ by computing the Serre spectral sequence of the Borel construction
\begin{align}
    C_2(T^d)\longrightarrow C_2(T^d)\times_{\Sigma_2}E\Sigma_2  \longrightarrow B\Sigma_2.
\end{align}
Since the action of $\Sigma_2$ on $C_2(T^d)$ is free, there is canonical homotopy equivalence $C_2(T^d)\times_{\Sigma_2}E\Sigma_2\rightarrow B_2(T^d)$. Hence, we obtain the homotopy fibration
\begin{align}\label{htpy_fib}
    C_2(T^d)\longrightarrow B_2(T^d)\longrightarrow B\Sigma_2.
\end{align}
Notice that the first map is the quotient map by the $\Sigma_2$ action, and the classifying space $B\Sigma_2$ is the projective space $\R P^\infty$. By computing the Serre spectral sequence of the homotopy fibration \ref{htpy_fib}, we obtain the following result.

\vspace{3mm}
\begin{theorem}\label{case_2_3}
We identify $\ccohomology{*}{\rp{\infty}}{\F_2}\cong \F_2[\alpha]$. When $d=2,3$, there are isomorphisms of $\F_2[\alpha]$-modules 
\begin{align}
    \ccohomology{*}{\uconfi{2}{\torus{2}}}{\F_2}\cong \F_2&[\alpha]/(\alpha^3)\oplus (\oplus_{j=1}^{2}\F_2 u_{1,j}) \oplus \F_2 u_{2,1} \oplus (\oplus_{j=1}^{2} \F_2[\alpha]/(\alpha^3) x_{2,j}),\\
    \ccohomology{*}{\uconfi{2}{\torus{3}}}{\F_2}\cong \F_2&[\alpha]/(\alpha^4) \oplus (\oplus_{j=1}^{3}\F_2 u_{1,j}) \oplus (\oplus_{j=1}^{6}\F_2 u_{2,j}) \oplus (\oplus_{j=1}^{3} \F_2[\alpha]/(\alpha^3) x_{2,j}) \\
    &\oplus (\oplus_{j=1}^{9} \F_2 u_{3,j}) \oplus (\oplus_{j=1}^{3} \F_2 u_{4,j})\oplus (\oplus_{j=1}^{3} \F_2[\alpha]/(\alpha^2) x_{4,j})
\end{align}
where $u_{i,j}$ and $x_{i,j}$ are elements with degree $i$.
\end{theorem}

\vspace{3mm}
To determine the module structure for the case of $d=3$, we use the following result:

\begin{theorem}\label{SWh}
    For any $d\geq 2$, the Stiefel-Whitney height of the $\Z_2$-bundle $\namelessmap{\confi{2}{\torus{d}}}{\uconfi{2}{\torus{d}}}$ is exactly $d$. Here the Stiefel-Whitney height of a $\Z_2$-bundle is the maximal integer $i$ such that $\alpha^i\neq 0$ where $\alpha$ is the first Stiefel-Whitney class of the bundle. 
\end{theorem}

The Stiefel-Whitney height of a manifold $M$ is an interesting invariant since it is related to the minimal dimension $n$ such that a triangulated skeleton of $M$ can be embedded in $\R^n$ (see \cite{kishimoto2023van}).

\vspace{3mm}
\section{Computation}

In this section, we compute the $\F_2[\sym{2}]$-module structure of $\ccohomology{i}{\confi{2}{\torus{d}}}{\F_2}$ and prove \cref{confi_rep} for $d\geq 1$. 

$\torus{d}=\R^d/\Z^d$ inherits the Lie group structure from $\R^d$. So we can define a map $\map{\varphi}{\torus{d}\times \torus{d}}{\torus{d}\times \torus{d}}$ by 
\begin{align}
    \varphi(x,y)=(x,y-x). 
\end{align} 
Note that the restriction $\map{\varphi}{\confi{2}{\torus{d}}}{\torus{d}\times (\torus{d}\setminus 0)}$ is a homeomorphism. Via this map, the desired representation is closely related to $\ccohomology{i}{\torus{d}\times \torus{d}}{\F_2}$. So we first compute the the $\F_2[\sym{2}]$-module structure of $\ccohomology{i}{\torus{d}\times \torus{d}}{\F_2}$.

\vspace{3mm}
\begin{lemma}\label{torus_rep}
    There is an isomorphism of $\F_2[\sym{2}]$-module
    \begin{align}
        \ccohomology{i}{\torus{d}\times \torus{d}}{\F_2} \cong 
        \begin{cases}
          {\F_2}^{\oplus \binom{d}{k}}\oplus \F_2[\sym{2}]^{\oplus \frac{1}{2}\left(\sum_{j=0}^{i}\binom{d}{j}\binom{d}{i-j}-\binom{d}{k}\right)} & \text{     if $i=2k$ for some $0\leq k\leq d$,}\\
          \F_2[\sym{2}]^{\oplus \frac{1}{2}\sum_{j=0}^{i}\binom{d}{j}\binom{d}{i-j} }& \text{     if $i=2k+1$ for some $0\leq k\leq d-1$,}\\
          0 & \text{  otherwise.}
        \end{cases}
    \end{align}
    Here we assume that $\binom{m}{n}=0$ when $m<n$.
    
\end{lemma}

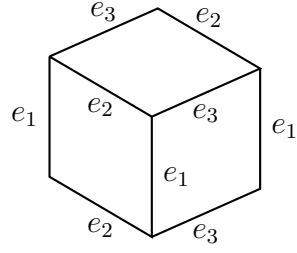
\begin{figure}
    \centering
    \begin{tikzpicture}[scale=0.8]
    \draw[thick]
    (0,0)-- node[left] (A) {$e_1$} (0,-2)--node[below] (B) {$e_2$}(1.7,-3)--node[right] (C) {$e_1$} (1.7,-1)--node[below] (D) {$e_2$}(0,0);
    \draw[thick]
    (0,0)--node[above] (E) {$e_3$}(1.8, 0.8)--node[above] (F) {$e_2$}(3.5, -0.2)--node[below] (G) {$e_3$}(1.7, -1);
    \draw[thick]
    (3.5,-0.2)--node[right] (H) {$e_1$}(3.5,-2.2)--node[below] (I) {$e_3$}(1.7,-3);
    \end{tikzpicture}
    \caption{Cell decomposition of $\torus{3}$}
    \label{cell_decomposition}
\end{figure}

\vspace{3mm}
Before the proof, we decompose $\torus{d}$ into cells, to use cellular homology and cohomology. We regard $\torus{d}$ to the quotient space obtained from the $d$-dimensional cube $[-\frac{1}{2}, \frac{1}{2}]^d$ by identifying each pair of opposite faces. Let $e_i$ denote the 1-cell of $\torus{d}$ corresponding to the edge from $(-\frac{1}{2},\dots, -\frac{1}{2})$ to $\underset{\hat{i}}{(-\frac{1}{2},\dots, -\frac{1}{2}, \frac{1}{2},-\frac{1}{2},\dots,-\frac{1}{2})}$(\cref{cell_decomposition}). Then each higher dimensional cell is a direct product of these 1-cells, and the only cell containing the origin is the $d$-cell $e_1\times \cdots \times e_d$. To simplify our notation, we omit the symbol $\times$ and the 0-cell for direct product. Note that $\chomology{i}{\torus{d}}{\Z}$ is a free abelian group generated by the $i$-dimensional cells $e_{k_1}e_{k_2}\cdots e_{k_i}$ for $0\leq k_1 < k_2 < \cdots <k_i\leq d$.

\begin{proof}[\textbf{Proof of \cref{torus_rep}}]
    By the universal coefficient theorem, there is a natural isomorphism
    \begin{align}
        \ccohomology{i}{\torus{d}}{\F_2}\cong \Hom{\Z}{\chomology{i}{\torus{d}}{\Z}}{\F_2}.
    \end{align}
    Let $e^{*}_{j}$ be the dual basis of $e_j$. By the Künneth theorem, we obtain a natural isomorphism
    \begin{align}
        \ccohomology{i}{\torus{d}\times \torus{d}}{\F_2}&\cong \bigoplus_{j=0}^{i}\ccohomology{j}{\torus{d}}{\F_2}\otimes \ccohomology{i-j}{\torus{d}}{\F_2}\\
        &\cong \bigoplus_{j=0}^{i}\Hom{\Z}{\chomology{j}{\torus{d}}{\Z}}{\F_2}\otimes \Hom{\Z}{\chomology{i-j}{\torus{d}}{\Z}}{\F_2}.
    \end{align}
    Via this isomorphism, the $\sym{2}$-action on $\ccohomology{i}{\torus{d}\times \torus{d}}{\F_2}$ is equivalent to swapping the left and right of the tensor product $\ccohomology{j}{\torus{d}}{\F_2}\otimes \ccohomology{i-j}{\torus{d}}{\F_2}$. Thus the subgroup generated by $e^{*}_{m_1}\cdots e^{*}_{m_s}\otimes e^{*}_{n_1}\cdots e^{*}_{n_t}$ and $e^{*}_{n_1}\cdots e^{*}_{n_t}\otimes e^{*}_{m_1}\cdots e^{*}_{m_s}$ is a subrepresentation of $\sym{2}$. Unless $e^{*}_{m_1}\cdots e^{*}_{m_s} = e^{*}_{n_1}\cdots e^{*}_{n_t}$, such representation is isomorphic to $\F_2[\sym{2}]=\{0,1,\sigma,1+\sigma\}$ by the homomorphism defined by 
    \begin{align}
        e^{*}_{m_1}\cdots e^{*}_{m_s}\otimes e^{*}_{n_1}\cdots e^{*}_{n_t}\mapsto 1, \quad e^{*}_{n_1}\cdots e^{*}_{n_t}\otimes e^{*}_{m_1}\cdots e^{*}_{m_s}\mapsto \sigma
    \end{align}
    When $e^{*}_{m_1}\cdots e^{*}_{m_s} = e^{*}_{n_1}\cdots e^{*}_{n_t}$, such a subrepresentation is trivial. 

    \vspace{3mm}
    If $i=2k+1$ for some $0\leq k\leq d$, there is no trivial subrepresentation. So we have an isomorphism of $\F_2[\sym{2}]$-module
    \begin{align}
        \ccohomology{i}{\torus{d}\times \torus{d}}{\F_2}\cong \F_2[\sym{2}]^{\oplus \frac{1}{2}\dim{\ccohomology{i}{\torus{d}\times \torus{d}}{\F_2}}}=\F_2[\sym{2}]^{\oplus \frac{1}{2}\sum_{j=0}^{i}\binom{d}{j}\binom{d}{i-j}}.
    \end{align}

    If $i=2k$ for some $0\leq k\leq d$, there are $\binom{d}{k}$ cells of the form $e_{m_1}\cdots e_{m_k}\otimes e_{m_1}\cdots e_{m_k}$, and for every other $i$-cell $e^{*}_{m_1}\cdots e^{*}_{m_s}\otimes e^{*}_{n_1}\cdots e^{*}_{n_{i-s}}$, there is unique swapped $i$-cell $e^{*}_{n_1}\cdots e^{*}_{n_{i-s}}\otimes e^{*}_{m_1}\cdots e^{*}_{m_s}$. So we have the direct sum decomposition of representation
    \begin{align}
        \ccohomology{i}{\torus{d}\times \torus{d}}{\F_2}\cong \F_2^{\oplus \binom{d}{k}}\oplus\F_2[\sym{2}]^{\oplus \frac{1}{2}\left(\sum_{j=0}^{i}\binom{d}{j}\binom{d}{i-j}-\binom{d}{k}\right)}.
    \end{align}

\end{proof}

We next show \cref{confi_rep}. The proof is supported by the commutative diagram
\begin{align}\label{ess_dia}
    \xymatrix@C=36pt{
    \confi{2}{\torus{d}} \ar[r]^-{\varphi}_-{\cong} \ar[d]_-{\iota} & \torus{d}\times (\torus{d}\setminus 0)\ar[d]^-{\iota'}\\
    \torus{d}\times \torus{d} \ar[r]_-{\varphi}^-{\cong}&\torus{d}\times \torus{d}
    }
\end{align}
where the vertical maps are the inclusions. 

\begin{proof}[\textbf{Proof of \cref{confi_rep}}]
    The diagram \ref{ess_dia} induces the commutative diagram of cohomology groups
    \begin{align}\label{cohom_dia}
        \xymatrix{
        \ccohomology{i}{\confi{2}{\torus{d}}}{\F_2}& \ccohomology{i}{\torus{d}\times (\torus{d}\setminus 0)}{\F_2}\ar[l]_-{\varphi^*} \\
        \ccohomology{i}{\torus{d}\times \torus{d}}{\F_2}\ar[u]^-{\iota^*}&\ccohomology{i}{\torus{d}\times \torus{d}}{\F_2}\ar[l]^-{\varphi^*}\ar[u]_-{\iota'^*}
        }
    \end{align}
    This diagram implies that the image of $e^{*}_{m_1}\cdots e^{*}_{m_s} \otimes e^{*}_{n_1}\cdots e^{*}_{n_t}\in \ccohomology{i}{\torus{d}\times \torus{d}}{\F_2}$ by $\iota^*$ generates $\ccohomology{i}{\confi{2}{\torus{d}}}{\F_2}$. 

    \vspace{3mm}
    If $i<d$, $\iota^*$ is an isomorphism since $\iota'^*$ and $\varphi^*$ are isomorphisms. So $\ccohomology{i}{\confi{2}{\torus{d}}}{\F_2}$ is isomorphic to $\ccohomology{i}{\torus{d}\times \torus{d}}{\F_2}$ as $\F_2$-module. Since the inclusion $\iota:\confi{2}{\torus{d}}\rightarrow \torus{d}\times \torus{d}$ is $\sym{2}$-equivariant, the diagram below is commutative, hence, $\ccohomology{i}{\confi{2}{\torus{d}}}{\F_2}$ is isomorphic to $\ccohomology{i}{\torus{d}\times \torus{d}}{\F_2}$ as $\F_2[\sym{2}]$-module ($\sigma$ is the generator of $\sym{2}$). 

    \begin{align}
        \xymatrix{
        \ccohomology{i}{\confi{2}{\torus{d}}}{\F_2}\ar[r]^-{\sigma}&\ccohomology{i}{\confi{2}{\torus{d}}}{\F_2}\\
        \ccohomology{i}{\torus{d}\times \torus{d}}{\F_2}\ar[r]_-{\sigma}\ar[u]^-{\iota^*}&\ccohomology{i}{\torus{d}\times \torus{d}}{\F_2}\ar[u]_-{\iota^*}
        }
    \end{align}

    If $i=d$, $\iota'^*:\ccohomology{i}{\torus{d}\times \torus{d}}{\F_2}\rightarrow \ccohomology{i}{\torus{d}\times (\torus{d}\setminus 0)}{\F_2}$ ,in \ref{cohom_dia}, is not an isomorphism, but it is a surjection whose kernel is the $1$-dimensional $\F_2$-vector space with basis $\varphi^*(1\otimes e^{*}_1\cdots e^{*}_d)$. The commutativity of \ref{cohom_dia} implies that $\map{\iota^*}{\ccohomology{i}{\torus{d}\times \torus{d}}{\F_2}}{\ccohomology{i}{\confi{2}{\torus{d}}}{\F_2}}$ is surjective. Hence we obtain an isomorphism of $\F_2[\sym{2}]$-module
    \begin{align}
        \ccohomology{i}{\confi{2}{\torus{d}}}{\F_2}\cong \ccohomology{i}{\torus{d}\times \torus{d}}{\F_2}/\ker{\iota^*}= \ccohomology{i}{\torus{d}\times \torus{d}}{\F_2}/\varphi^*(\langle 1\otimes e^{*}_{1}\cdots e^{*}_{d}\rangle).
    \end{align}  
    Here we note that $\ker{\iota^*}$ is a $\F_2[\sym{2}]$-submodule of $\ccohomology{i}{\torus{d}\times \torus{d}}{\F_2}$ since $\iota^*$ is $\sym{2}$-equivariant. To determine the $\F_2[\sym{2}]$-module structure, we next compute the induced map $\varphi^*:\ccohomology{i}{\torus{d}\times \torus{d}}{\F_2}\rightarrow \ccohomology{i}{\torus{d}\times \torus{d}}{\F_2}$.

    \vspace{3mm}
    Since all the cells of dimension greater than 1 of $\torus{d}\times \torus{d}$ are expressed as a product cell of 1-cells, so all the cohomology classes of such dimension of $\torus{d}\times \torus{d}$ are generated by cup products of some classes of 1-cocycles. Since $\varphi^*:\ccohomology{*}{\torus{d}\times \torus{d}}{\F_2}\rightarrow \ccohomology{*}{\torus{d}\times \torus{d}}{\F_2}$ is a ring homomorphism, $\varphi^*$ is completely determined by the restriction on the first cohomology group $\ccohomology{1}{\torus{d}\times \torus{d}}{\F_2}$. By the universal coefficient theorem, it is the dual map $(\varphi_*)_\#$ of $\varphi_*:\chomology{1}{\torus{d}\times \torus{d}}{\Z}\rightarrow \chomology{1}{\torus{d}\times \torus{d}}{\Z}$. 
    \begin{align}
        \xymatrix{
        \ccohomology{1}{\torus{d}\times \torus{d}}{\F_2} \ar[r]^-{\cong} \ar[d]_-{\varphi^*}& \Hom{\Z}{\chomology{1}{\torus{d}\times \torus{d}}{\Z}}{\F_2}\ar[d]^-{(\varphi_*)^\#}\\
        \ccohomology{1}{\torus{d}\times \torus{d}}{\F_2} \ar[r]_-{\cong} & \Hom{\Z}{\chomology{1}{\torus{d}\times \torus{d}}{\Z}}{\F_2}
        }
    \end{align}

    By the definition of $\varphi:\torus{d}\times \torus{d}\rightarrow \torus{d}\times \torus{d}$, we have
    \begin{align}
        \varphi_*(e_i\otimes *)&= (e_i \otimes *) - (*\otimes e_i),\\
        \varphi_*(*\otimes e_i)&= *\otimes e_i,
    \end{align}
    where $*$ is the 0-cell of $\torus{d}$. So the dual map is determined by
    \begin{align}
        \varphi^*(e_{i}^{*}\otimes 1)&= e_{i}^{*}\otimes 1,\\
        \varphi^*(1\otimes e_{i}^{*})&= (e_{i}^{*}\otimes 1)- (1\otimes e_{i}^{*})=(e_{i}^{*}\otimes 1)+ (1\otimes e_{i}^{*}).
    \end{align}

    \vspace{3mm}
    We return to the consideration for the case $i=d$. By the preceding calculation, we have
    \begin{align}\label{relation}
        \varphi^*(1\otimes e^{*}_{1}\cdots e^{*}_{d})&=\prod_{j=1}^{d}\varphi^*(1\otimes e_{j}^{*})\\
        &=\prod_{j=1}^{d}(e_{j}^{*}\otimes 1 + 1\otimes e_{j}^{*})\\
        &=\sum_{l=0}^{d}\sum_{1\leq j_{k_1}< \cdots < j_{k_l}\leq d} e_{1}^{*} \cdots \hat{e_{j_{k_1}}^{*}}\cdots \hat{e_{j_{k_l}}^{*}}\cdots e_{d}^{*} \otimes e_{j_{k_1}}^{*}\cdots e_{j_{k_l}}^{*}.
    \end{align}
    Thus $\ccohomology{i}{\torus{d}\times \torus{d}}{\F_2}/\varphi^*(\langle 1\otimes e^{*}_{1}\cdots e^{*}_{d}\rangle)$ is the $\F_2[\sym{2}]$-module generated by $e^{*}_{m_1}\cdots e^{*}_{m_s} \otimes e^{*}_{n_1}\cdots e^{*}_{n_{d-s}}$ with one relation
    \begin{align}
        \sum_{l=0}^{d}\sum_{1\leq j_{k_1}< \cdots < j_{k_l}\leq d} e_{1}^{*} \cdots \hat{e_{j_{k_1}}^{*}}\cdots \hat{e_{j_{k_l}}^{*}}\cdots e_{d}^{*} \otimes e_{j_{k_1}}^{*}\cdots e_{j_{k_l}}^{*}=0.
    \end{align}
    Notice that every summand $e_{1}^{*} \cdots \hat{e_{j_{k_1}}^{*}}\cdots \hat{e_{j_{k_l}}^{*}}\cdots e_{d}^{*} \otimes e_{j_{k_1}}^{*}\cdots e_{j_{k_l}}^{*}$ corresponds to the unique choice of the indices in the left side of the tensor product, which is equivalent to the choice of the indices in the right side, and this is one-to-one. Thus choosing a half of the summands, we obtain a trivial suburepresentation of $\sym{2}$. More explicitly, for example, define $x$ by the formula
    \begin{align}
        x=\begin{dcases}
            \sum_{l=0}^{(d/2)-1} \sum_{1\leq j_{k_1}< \cdots < j_{k_l}\leq d} e_{1}^{*} \cdots \hat{e_{j_{k_1}}^{*}}\cdots \hat{e_{j_{k_l}}^{*}}\cdots e_{d}^{*} \otimes e_{j_{k_1}}^{*}\cdots e_{j_{k_l}}^{*}\\
            \quad\quad\quad +\sum_{1\leq j_{k_1}< \cdots < j_{k_{d/2}}< d} e_{1}^{*} \cdots \hat{e_{j_{k_1}}^{*}}\cdots \hat{e_{j_{k_{d/2}}}^{*}}\cdots e_{d}^{*} \otimes e_{j_{k_1}}^{*}\cdots e_{j_{k_{d/2}}}^{*} & \text{  if $d$ is even,}\\
            \sum_{l=0}^{(d-1)/2} \sum_{1\leq j_{k_1}< \cdots < j_{k_l}\leq d} e_{1}^{*} \cdots \hat{e_{j_{k_1}}^{*}}\cdots \hat{e_{j_{k_l}}^{*}}\cdots e_{d}^{*} \otimes e_{j_{k_1}}^{*}\cdots e_{j_{k_l}}^{*}&\text{  if $d$ is odd.}
        \end{dcases}
    \end{align}
    Then $x$ is a nontrivial element satisfying $\sigma(x)+x=0$, which is equivalent to that $\sigma(x)=x$. So we have that, in $\ccohomology{i}{\torus{d}\times \torus{d}}{\F_2}/\varphi^*(\langle 1\otimes e^{*}_{1}\cdots e^{*}_{d}\rangle)$, the regular representation $\langle e^{*}_{1}\cdots e^{*}_{d}\otimes 1, 1\otimes e^{*}_{1}\cdots e^{*}_{d}\rangle$ in $\ccohomology{i}{\torus{d}\times \torus{d}}{\F_2}$ has changed into the trivial representation $\langle x\rangle$, and all the other subrepresentations of $\ccohomology{i}{\torus{d}\times \torus{d}}{\F_2}$ are the same. Therefore, by \cref{torus_rep}, we have
    \begin{align}
        \ccohomology{d}{\confi{2}{\torus{d}}}{\F_2}&\cong \ccohomology{d}{\torus{d}\times \torus{d}}{\F_2}/\varphi^*(\langle 1\otimes e^{*}_{1}\cdots e^{*}_{d}\rangle)\\
        &\cong 
        \begin{cases}
          {\F_2}^{\oplus \binom{d}{k}+1}\oplus \F_2[\sym{2}]^{\oplus \frac{1}{2}\left(\sum_{j=0}^{i}\binom{d}{j}\binom{d}{d-j}-\binom{d}{k}\right)-1} & \text{     if $d=2k$ for some $k\geq 0$,}\\
          \F_2\oplus\F_2[\sym{2}]^{\oplus \frac{1}{2}\sum_{j=0}^{i}\binom{d}{j}\binom{d}{d-j} -1}& \text{     if $d=2k+1$ for some $k\geq 0$}.
        \end{cases}
    \end{align}

    It remains the case $i>d$. By the same diagram \ref{cohom_dia}, $\ccohomology{i}{\confi{2}{\torus{d}}}{\F_2}$ is isomorphic to $\ccohomology{i}{\torus{d}\times \torus{d}}{\F_2}/\ker{\iota^*}$. Now the kernel of $\iota^*$ is generated by the images
    \begin{align}
        \varphi^*(e^{*}_{m_1}\cdots e^{*}_{m_{i-d}} \otimes e^{*}_{1}\cdots e^{*}_{d})&=\varphi^*(e^{*}_{m_1}\cdots e^{*}_{m_{i-d}} \otimes 1) \varphi^*(1\otimes e^{*}_{1}\cdots e^{*}_{d})\\
        &=\left(\prod_{j=1}^{i-d} \varphi^*(e^{*}_{m_{j}}\otimes 1) \right)\varphi^*(1\otimes e^{*}_{1}\cdots e^{*}_{d})\\
        &=\left(\prod_{j=1}^{i-d}(e^{*}_{m_{j}}\otimes 1) \right)\varphi^*(1\otimes e^{*}_{1}\cdots e^{*}_{d})\\
        &=(e^{*}_{m_1}\cdots e^{*}_{m_{i-d}} \otimes 1)\varphi^*(1\otimes e^{*}_{1}\cdots e^{*}_{d}).
    \end{align}
    By \cref{relation}, in $\ccohomology{i}{\torus{d}\times \torus{d}}{\F_2}/\ker{\iota^*}$, there are $\binom{d}{i-d}$ relations, each of which is given by
    \begin{align}
        (e^{*}_{m_1}\cdots e^{*}_{m_{i-d}} \otimes 1)\left(\sum_{l=0}^{d}\sum_{1\leq j_{k_1}< \cdots < j_{k_l}\leq d} e_{1}^{*} \cdots \hat{e_{j_{k_1}}^{*}}\cdots \hat{e_{j_{k_l}}^{*}}\cdots e_{d}^{*} \otimes e_{j_{k_1}}^{*}\cdots e_{j_{k_l}}^{*}\right)=0.
    \end{align}
    Since $(e_{j}^{*})^2=0$ foe each $j$, this relation is equivalent to
    \begin{align}
        \left(\sum_{} e_{1}^{*} \cdots \hat{e_{j_{k_1}}^{*}}\cdots \hat{e_{j_{k_l}}^{*}}\cdots e_{d}^{*} \otimes e_{j_{k_1}}^{*}\cdots \hat{e^{*}_{m_1}}\cdots \hat{e^{*}_{m_{i-d}}}\cdots e_{j_{k_l}}^{*}\right)(e^{*}_{m_1}\cdots e^{*}_{m_{i-d}}\otimes e^{*}_{m_1}\cdots e^{*}_{m_{i-d}})=0
    \end{align}
    where $l$ ranges $0$ to $d$, and $j_{k_1}, \dots,j_{k_l}$ range $1\leq j_{k_1}< \cdots < j_{k_l}\leq d$ under $j_{k_s}\neq m_1,\dots,m_{i-d}$ for each $1\leq s\leq l$. Notice that, in the summation, each summand does not contain $e^{*}_{m_1},\dots , e^{*}_{m_{i-d}}$, and the indices in each summand do not duplicate. Since the summation is summing up all of such summands, we can pair the summands of the form $e^{*}_{m_1}\cdots e^{*}_{m_s}\otimes e^{*}_{n_1}\cdots e^{*}_{n_t}$ and $e^{*}_{n_1}\cdots e^{*}_{n_t}\otimes e^{*}_{m_1}\cdots e^{*}_{m_s}$, then we can construct a nontrivial element $x'$ so that $\sigma(x')=x'$ as well as the previous $x$. In the end, in $\ccohomology{i}{\torus{d}\times \torus{d}}{\F_2}/\ker{\iota^*}$,  $\binom{d}{i-d}$ regular representations $\langle e^{*}_{1}\cdots e^{*}_{d}\otimes e^{*}_{m_1}\cdots e^{*}_{m_{i-d}} , e^{*}_{m_1}\cdots e^{*}_{m_{i-d}} \otimes e^{*}_{1}\cdots e^{*}_{d}\rangle$ in $\ccohomology{i}{\torus{d}\times \torus{d}}{\F_2}$ have changed into trivial representations. However, when $i=2d$, the relation means that $\ker{\iota^*}= \ccohomology{2d}{\torus{d}\times \torus{d}}{\F_2}$. Therefore, by \cref{torus_rep}, we finally obtain 
    \begin{align}
        &\ccohomology{i}{\confi{2}{\torus{d}}}{\F_2} \cong \ccohomology{i}{\torus{d}\times \torus{d}}{\F_2}/\ker{\iota^*}\\
        &\cong\begin{cases}
          {\F_2}^{\oplus \binom{d}{k}+\binom{d}{i-d}}\oplus \F_2[\sym{2}]^{\oplus \frac{1}{2}\left(\sum_{j=0}^{i}\binom{d}{j}\binom{d}{i-j}-\binom{d}{k}\right)-\binom{d}{i-d}} & \text{     if $d\leq i=2k<d$,}\\
          \F_2^{\oplus \binom{d}{i-d}}\F_2[\sym{2}]^{\oplus \frac{1}{2}\left(\sum_{j=0}^{i}\binom{d}{j}\binom{d}{i-j} \right) - \binom{d}{i-d}}& \text{     if $d\leq i=2k+1<d$,}\\
          0& \text{    if $i\geq 2d$.}
        \end{cases}
    \end{align}
    
\end{proof}

\vspace{5mm}
In particular, when $d\leq 3$, we have 
\begin{align}
    \rccohomology{i}{\confi{2}{\torus{1}}}{\F_2}=\rccohomology{i}{\confi{2}{S^1}}{\F_2}&\cong 
    \begin{cases}
         0&\text{$i=0$,}\\
        \F_2[\sym{2}] &\text{$i=1$,}\\
        0&\text{otherwise.}\\
    \end{cases}\\
    \rccohomology{i}{\confi{2}{\torus{2}}}{\F_2}&\cong 
    \begin{cases}
         0&\text{$i=0$,}\\
        \F_2[\sym{2}]^{\oplus 2} &\text{$i=1$,}\\
        \F_2^{\oplus 3} \oplus \F_2[\sym{2}] &\text{$i=2$,}\\
        \F_2^{\oplus 2} &\text{$i=3$,}\\
        0&\text{otherwise.}\\
    \end{cases}\\
    \rccohomology{i}{\confi{2}{\torus{3}}}{\F_2}&\cong
    \begin{cases}
        0 &\text{$i=0$,}\\
        \F_2[\sym{2}]^{\oplus 3} &\text{$i=1$,}\\
        \F_2^{\oplus 3} \oplus \F_2[\sym{2}]^{\oplus 6} &\text{$i=2$,}\\
        \F_2^{\oplus 1} \oplus \F_2[\sym{2}]^{\oplus 9} &\text{$i=3$,}\\
        \F_2^{\oplus 6} \oplus \F_2[\sym{2}]^{\oplus 3} &\text{$i=4$,}\\
        \F_2^{\oplus 3} &\text{$i=5$,}\\
        0&\text{otherwise.}
    \end{cases}
\end{align}
for the reduced cohomology groups. In these cases, the rank of the trivial representation in $H^i$ and that of the regular one in $H^{d-i}$ are symmetrically swapped. Although such `duality'-like situation does not occur for $d\geq 4$, this is a little interesting.

\vspace{3mm}
\section{Applications}

In this section, we prove \cref{case_2_3,SWh}. We first see the $\ccohomology{*}{\rp{\infty}}{\F_2}$-module structure of $\ccohomology{i}{\uconfi{2}{\torus{2}}}{\F_2}$ since it can be determined directly. The $E_2$-term of the Serre spectral sequence for the homotopy fibration \ref{htpy_fib} is given by 
\begin{align}
    E_{2}^{p,q}=\ccohomology{p}{B\sym{2}}{\lcohomology{q}{\confi{2}{\torus{2}}}{\F_2}}
\end{align}
where $\lcohomology{q}{\confi{2}{\torus{d}}}{\F_2}$ denotes the local system over $\rp{\infty}$ associated to the homotopy fibration \ref{htpy_fib}. By the preceding result, we have
\begin{align}
    E_{2}^{p,q}=
    \begin{cases}
        \ccohomology{p}{\R P^\infty}{\F_2}& \text{ $q=0$,}\\
        \ccohomology{p}{\R P^\infty}{\F_2[\sym{2}]}^{\oplus 2}& \text{ $q=1$,}\\
        \ccohomology{p}{\R P^\infty}{\F_2}^{\oplus 3}\oplus \ccohomology{p}{\R P^\infty}{\F_2[\sym{2}]}^{\oplus 2}& \text{ $q=2$,}\\
        \ccohomology{p}{\R P^\infty}{\F_2}^{\oplus 2}& \text{ $q=3$,}\\
        0& \text{ otherwise.}
    \end{cases}
\end{align}
In general, for a connected CW complex $X$ with finite fundamental group $\pi$ and the universal cover $\tilde{X}$, there is an isomorphism $\ccohomology{i}{\tilde{X}}{R}\cong \ccohomology{i}{X}{R[\pi]}$ (see the proof of Proposition 3H.5 in \cite{MR1867354}). Since $\abhtpygrp{1}{\rp{\infty}}\cong \sym{2}$, we have $\ccohomology{i}{\rp{\infty}}{\F_2[\sym{2}]}\cong \ccohomology{i}{\sphere{\infty}}{\F_2}$ by considering the universal cover $\namelessmap{\sphere{\infty}}{\rp{\infty}}$. This isomorphism is derived from the isomorphism of cochain complexes $\namelessmap{\cellcochain{*}{\sphere{\infty}}{\F_2}}{\Hom{\F_2[\sym{2}]}{\cellchain{*}{\sphere{\infty}}}{\F_2[\sym{2}]}}$ where $\cellchain{*}{\sphere{\infty}}$ and $\cellcochain{*}{\sphere{\infty}}{\F_2}$ denote the cellular chain and cochain complexes, with coefficients in integer and $\F_2$, respectively. More explicitly, this map is given by
\begin{align}
\varphi\mapsto\left(\hat{\varphi}\colon x\mapsto\sum_{\gamma\in\pi_1}\varphi(\gamma^{-1}x)\gamma\right)
\end{align}
where $\pi_1$ is the fundamental group of $\rp{\infty}$. Hence, 
\begin{align}
    \ccohomology{p}{\rp{\infty}}{\F_2[\sym{2}]}\cong \ccohomology{p}{\sphere{\infty}}{\F_2}\cong
    \begin{cases}
        \F_2 & \text{if $p=0$,}\\
        0 & \text{otherwise.}
    \end{cases}
\end{align}
Since the generator of $\ccohomology{0}{\sphere{\infty}}{\F_2}$ is the $0$-cochain $1$ of $\sphere{\infty}$ corresponding to $\hat{1}\colon e^0\mapsto 1+\sigma$, thus $\ccohomology{0}{\rp{\infty}}{\F_2[\sym{2}]}$ is generated by $e_{i}^{*}\otimes 1+1\otimes e_{i}^{*}$. Thus the $E_2$-term is described as follows:
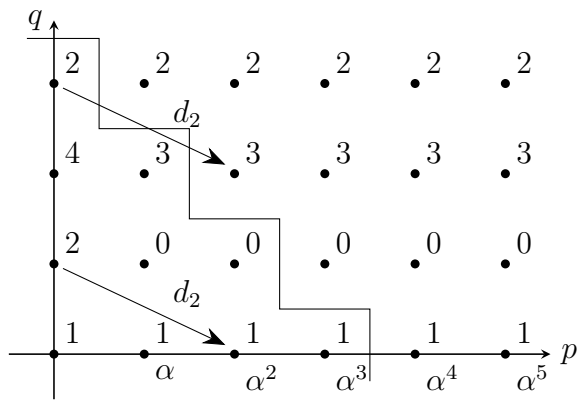
\begin{figure}[H]
    \centering
    \begin{tikzpicture}[xscale=1.2,yscale=1.2]
        \coordinate(O)at(0,0);
        \coordinate(XS)at(-0.5,0);
        \coordinate(XL)at(5.5,0);
        \coordinate(YS)at(0,-0.5);
        \coordinate(YL)at(0,3.7);
        \draw[semithick,->,>=stealth](XS)--(XL)node[right]{$p$};
        \draw[semithick,->,>=stealth](YS)--(YL)node[left]{$q$};
        \coordinate(P)at(1,0);
        \coordinate(Q)at(0,1);

        \foreach\k in{0,1,2,3,4,5}\foreach\l in{0,1,2,3}\fill($(O)+\k*(P)+\l*(Q)$)circle(0.05);
        \foreach\k in{0,1,2,3,4,5}
        \coordinate[label=above right:$1$](X\k)at($\k*(P)$);
        \coordinate[label=below right:$\alpha$](X)at(1,0);
        \foreach\k in{2,3,4,5}
        \coordinate[label=below right:$\alpha^\k$](X\k)at($\k*(P)$);
        \coordinate[label=above right:$2$](Y)at(0,1);
        \foreach\k in{1,2,3,4,5}
        \coordinate[label=above right:$0$](Y\k)at($\k*(P)+(Q)$);
        \coordinate[label=above right:$4$](Z)at(0,2);
        \foreach\k in{1,2,3,4,5}\coordinate[label=above right:$3$](Z\k)at($\k*(P)+(Z)$);
        \coordinate(W)at(0,3);
        \foreach\k in{0,1,2,3,4,5}\coordinate[label=above right:$2$](W\k)at($\k*(P)+(W)$);

        \draw(-0.3,3.5)--(0.5,3.5)--(0.5,2.5)--(1.5,2.5)--(1.5,1.5)--(2.5,1.5)--(2.5,0.5)--(3.5,0.5)--(3.5,-0.3);

        \draw [-{Stealth[length=3mm]}] (0.1,0.95) -- (1.9,0.1);
        \coordinate[label=above right:$d_2$](d)at(1.2,0.4);
        \draw [-{Stealth[length=3mm]}] (0.1,2.95) -- (1.9,2.1);
        \coordinate[label=above right:$d_2$](d)at(1.2,2.4);

    \end{tikzpicture}
    \caption{$E_2$-term}
    
\end{figure}

\noindent
Here $\alpha$ denotes the generator of the polynomial ring $\ccohomology{*}{\rp{\infty}}{\F_2}$. Since $\uconfi{2}{\torus{2}}$ is a $4$-dimensional open manifold, the region of total degree $4$ or greater must vanish in $E_\infty$-term. So the upper $d_2$ in the figure is injective. On the other hand, the lower $d_2$ is zero because if $d_2(u)=\alpha^2$, then $\alpha^3=d_2(u)\alpha=d_2(u\alpha)=d_2(0)=0$. This contradicts to $\alpha^3\neq0$. Thus the $E_3$-term is given by
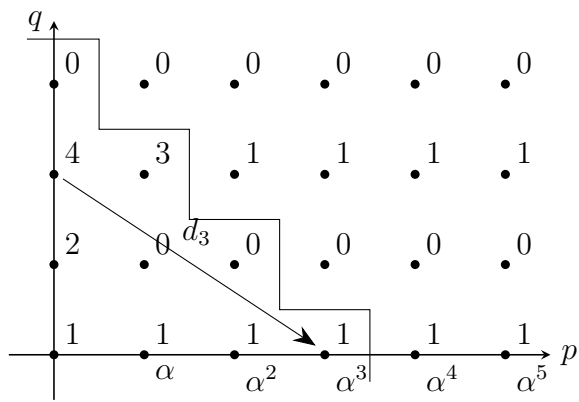
\begin{figure}[H]
    \centering
    \begin{tikzpicture}[xscale=1.2,yscale=1.2]
        \coordinate(O)at(0,0);
        \coordinate(XS)at(-0.5,0);
        \coordinate(XL)at(5.5,0);
        \coordinate(YS)at(0,-0.5);
        \coordinate(YL)at(0,3.7);
        \draw[semithick,->,>=stealth](XS)--(XL)node[right]{$p$};
        \draw[semithick,->,>=stealth](YS)--(YL)node[left]{$q$};
        \coordinate(P)at(1,0);
        \coordinate(Q)at(0,1);

        \foreach\k in{0,1,2,3,4,5}\foreach\l in{0,1,2,3}\fill($(O)+\k*(P)+\l*(Q)$)circle(0.05);
        \foreach\k in{0,1,2,3,4,5}
        \coordinate[label=above right:$1$](X\k)at($\k*(P)$);
        \coordinate[label=below right:$\alpha$](X)at(1,0);
        \foreach\k in{2,3,4,5}
        \coordinate[label=below right:$\alpha^\k$](X\k)at($\k*(P)$);
        \coordinate[label=above right:$2$](Y)at(0,1);
        \foreach\k in{1,2,3,4,5}
        \coordinate[label=above right:$0$](Y\k)at($\k*(P)+(Q)$);
        \coordinate[label=above right:$4$](Z)at(0,2);
        \coordinate[label=above right:$3$](Z1)at(1,2);
        \foreach\k in{2,3,4,5}\coordinate[label=above right:$1$](Z\k)at($\k*(P)+(Z)$);
        \coordinate(W)at(0,3);
        \foreach\k in{0,1,2,3,4,5}\coordinate[label=above right:$0$](W\k)at($\k*(P)+(W)$);

        \draw(-0.3,3.5)--(0.5,3.5)--(0.5,2.5)--(1.5,2.5)--(1.5,1.5)--(2.5,1.5)--(2.5,0.5)--(3.5,0.5)--(3.5,-0.3);

        \draw [-{Stealth[length=3mm]}] (0.1,1.95) -- (2.9,0.1);
        \coordinate[label=above right:$d_3$](d)at(1.3,1.1);

    \end{tikzpicture}
    \caption{$E_3$-term}
\end{figure}
\noindent
If the $d_3$ in the figure is zero, then $\alpha^4$ remains in $E_\infty$. Thus $d_3$ is nonzero, and we have the $E_\infty$ term indicated in the following figure:
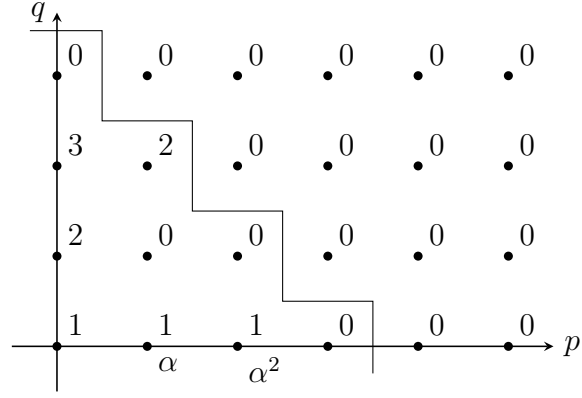
\begin{figure}[H]
    \centering
    \begin{tikzpicture}[xscale=1.2,yscale=1.2]
        \coordinate(O)at(0,0);
        \coordinate(XS)at(-0.5,0);
        \coordinate(XL)at(5.5,0);
        \coordinate(YS)at(0,-0.5);
        \coordinate(YL)at(0,3.7);
        \draw[semithick,->,>=stealth](XS)--(XL)node[right]{$p$};
        \draw[semithick,->,>=stealth](YS)--(YL)node[left]{$q$};
        \coordinate(P)at(1,0);
        \coordinate(Q)at(0,1);

        \foreach\k in{0,1,2,3,4,5}\foreach\l in{0,1,2,3}\fill($(O)+\k*(P)+\l*(Q)$)circle(0.05);
        \foreach\k in{0,1,2}
        \coordinate[label=above right:$1$](X\k)at($\k*(P)$);
        \foreach\k in{3,4,5}
        \coordinate[label=above right:$0$](X\k)at($\k*(P)$);
        \coordinate[label=below right:$\alpha$](X)at(1,0);
        \foreach\k in{2}
        \coordinate[label=below right:$\alpha^\k$](X\k)at($\k*(P)$);
        \coordinate[label=above right:$2$](Y)at(0,1);
        \foreach\k in{1,2,3,4,5}
        \coordinate[label=above right:$0$](Y\k)at($\k*(P)+(Q)$);
        \coordinate[label=above right:$3$](Z)at(0,2);
        \coordinate[label=above right:$2$](Z1)at(1,2);
        \foreach\k in{2,3,4,5}\coordinate[label=above right:$0$](Z\k)at($\k*(P)+(Z)$);
        \coordinate(W)at(0,3);
        \foreach\k in{0,1,2,3,4,5}\coordinate[label=above right:$0$](W\k)at($\k*(P)+(W)$);

        \draw(-0.3,3.5)--(0.5,3.5)--(0.5,2.5)--(1.5,2.5)--(1.5,1.5)--(2.5,1.5)--(2.5,0.5)--(3.5,0.5)--(3.5,-0.3);

    \end{tikzpicture}
    \caption{$E_\infty$-term}
\end{figure}
\noindent
So we can see that the Stiefel-Whitney height of the bundle $\namelessmap{\confi{2}{\torus{2}}}{\uconfi{2}{\torus{2}}}$ is exactly $2$. 

\vspace{5mm}
Before the computation for the case of $d=3$, we prove \cref{SWh}. This follows from the observation of the Serre spectral sequence and the fact that $\torus{d}$ can be embedded into $\R^{d+1}$.

\begin{proof}[\textbf{Proof of \cref{SWh}}]
    An embedding $\namelessmap{\torus{d}}{\R^{d+1}}$ induces a $\sym{2}$-equivariant embedding $\namelessmap{\confi{2}{\torus{d}}}{\confi{2}{\R^{d+1}}}$. Since there is a $\sym{2}$-equivariant homotopy equivalence $\confi{2}{\R^{d+1}}\simeq \sphere{d}$, we have $\uconfi{2}{\R^{d+1}}\simeq \rp{d}$. Thus the Stiefel-Whitney height of $\uconfi{2}{\torus{d}}$ is $d$ or less. On the other hand, the inclusion $\namelessmap{\confi{2}{\torus{d}}}{\torus{d}\times \torus{d}}$ induces a bundle map between the Borel constructions. We denote the Serre spectral sequence of the Borel construction 
    \begin{align}\label{borel2}
        \xymatrix{
            \torus{d}\times \torus{d}\ar[r]&(\torus{d}\times \torus{d})\times_{\sym{2}}E\sym{2}\ar[r]&B\sym{2}
        }
    \end{align}
    by $\{E',d'\}$. So its second page is indicated as
    \begin{align}
        E^{\prime p,q}_{2}=\ccohomology{p}{B\sym{2}}{\lcohomology{q}{\torus{d}\times \torus{d}}{\F_2}}.
    \end{align}
    Then, by \cref{confi_rep}, the induced homomorphism $\namelessmap{E'_r}{E_r}$ gives surjections restricted on the part of $q\leq d$. In fact, every differential of $E'$ whose source is in such part is zero. Thus that of $E$ so is. Therefore the Stiefel-Whitney classes $\alpha^i$ remain for $i=0,1,\dots,d$, and the Stiefel-Whitney height of $\uconfi{2}{\torus{d}}$ is $d$ or greater.

    \vspace{3mm}
    To complete the proof, we show that the Serre spectral sequence $\{E',d'\}$ associated to \ref{borel2} collapses at $E'_2$. This is checked by cup product of $E'$. First, we notice that $d'_2(u)=0$ if $u$ is in a regular representation part $\ccohomology{p}{\rp{\infty}}{\F_2[\sym{2}]}$ of $E_{2}^{\prime p,q}$. Indeed, if $d'(u)\neq 0$, we have $d'(u)\alpha=d'(u\alpha)\neq 0$. However, $u\alpha$ is zero because $u\alpha$ is an element in the regular representation part $\ccohomology{p+1}{\rp{\infty}}{\F_2[\sym{2}]\otimes_{\F_2} \F_2}$ of $E_{2}^{\prime p+1,q}$. Thus we have $d'_2=0$.
    
    On the other hand, the basis of the trivial representation part of $\ccohomology{*}{B\sym{2}}{\lcohomology{q}{\torus{d}\times \torus{d}}{\F_2}}$ consists of $e^{*}_{i}\otimes e^{*}_{i},\alpha$ and all of their products. However, we have $d'_3(e^{*}_{i}\otimes e^{*}_{i})=0$ as follows. Since the $\sym{2}$-action on $\torus{d}\times \torus{d}$ has a fixed point, a section of the Borel construction \ref{borel2} is obtained as the middle vertical map in the diagram
    \begin{align}
        \xymatrix@R=12pt{
        \torus{d}\times \torus{d} \ar[r] & (\torus{d}\times \torus{d})\times_{\sym{2}} E\sym{2} \ar[r] & B\sym{2}\\
        \ast \ar[r] \ar[u] & \ast \times_{\sym{2}} E\sym{2} = B\sym{2} \ar@{=}[r] \ar[u] & B\sym{2} \ar@{=}[u]
        }
    \end{align}
    So any power $\alpha^i$ of the first Stiefel-Whitney class $\alpha$ cannot be zero in $\ccohomology{*}{(\torus{d}\times \torus{d})\times_{\sym{2}} E\sym{2}}{\F_2}$. Hence, $d'_3(e^{*}_{i}\otimes e^{*}_{i})=0$, and the Leibniz rule implies that $d'_3=0$. Inductively, we have $d'_r=0$ for arbitrary $r\geq 2$.
\end{proof}

Now we can determine the $E_\infty$-term of the Serre spectral sequence for $d=3$. The $E_2$-term is indicated as the figure below:
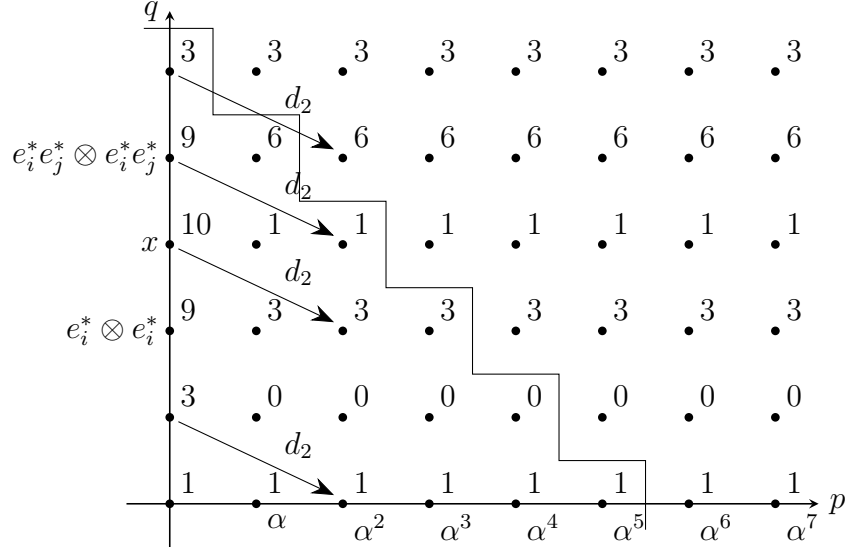
\begin{figure}[H]
    \centering
    \begin{tikzpicture}[xscale=1.15,yscale=1.15]
        \coordinate(O)at(0,0);
        \coordinate(XS)at(-0.5,0);
        \coordinate(XL)at(7.5,0);
        \coordinate(YS)at(0,-0.5);
        \coordinate(YL)at(0,5.7);
        \draw[semithick,->,>=stealth](XS)--(XL)node[right]{$p$};
        \draw[semithick,->,>=stealth](YS)--(YL)node[left]{$q$};
        \coordinate(P)at(1,0);
        \coordinate(Q)at(0,1);

        \foreach\k in{0,1,2,3,4,5,6,7}\foreach\l in{0,1,2,3,4,5}\fill($(O)+\k*(P)+\l*(Q)$)circle(0.05);
        
        \foreach\k in{0,1,2,3,4,5,6,7}
        \coordinate[label=above right:$1$](X\k)at($\k*(P)$);
        
        \coordinate[label=below right:$\alpha$](X)at(1,0);
        \foreach\k in{2,3,4,5,6,7}
        \coordinate[label=below right:$\alpha^\k$](X\k)at($\k*(P)$);
        
        \coordinate[label=above right:$3$](Y)at(0,1);
        
        \foreach\k in{1,2,3,4,5,6,7}
        \coordinate[label=above right:$0$](Y\k)at($\k*(P)+(Q)$);
        
        \coordinate[label=above right:$9$, label=left:$e^{*}_{i}\otimes e^{*}_{i}$](Z)at(0,2);
        \foreach\k in{1,2,3,4,5,6,7}\coordinate[label=above right:$3$](Z\k)at($\k*(P)+(Z)$);
        
        \coordinate[label=above right:$10$, label=left:$x$](W)at(0,3);
        
        \foreach\k in{1,2,3,4,5,6,7}\coordinate[label=above right:$1$](W\k)at($\k*(P)+(W)$);

        \coordinate[label=above right:$9$, label=left:$e^{*}_{i}e^{*}_{j}\otimes e^{*}_{i}e^{*}_{j}$](R)at(0,4);
        
        \foreach\k in{1,2,3,4,5,6,7}\coordinate[label=above right:$6$](R\k)at($\k*(P)+(0,4)$);

        \foreach\k in{0,1,2,3,4,5,6,7}\coordinate[label=above right:$3$](R\k)at($\k*(P)+(0,5)$);

        \draw(-0.3,5.5)--(0.5,5.5)--(0.5,4.5)--(1.5,4.5)--(1.5,3.5)--(2.5,3.5)--(2.5,2.5)--(3.5,2.5)--(3.5,1.5)--(4.5,1.5)--(4.5,0.5)--(5.5,0.5)--(5.5,-0.3);

        \draw [-{Stealth[length=3mm]}] (0.1,0.95) -- (1.9,0.1);
        \coordinate[label=above right:$d_2$](d)at(1.2,0.4);
        \draw [-{Stealth[length=3mm]}] (0.1,2.95) -- (1.9,2.1);
        \coordinate[label=above right:$d_2$](d)at(1.2,2.4);
        \draw [-{Stealth[length=3mm]}] (0.1,3.95) -- (1.9,3.1);
        \coordinate[label=above right:$d_2$](d)at(1.2,3.4);
        \draw [-{Stealth[length=3mm]}] (0.1,4.95) -- (1.9,4.1);
        \coordinate[label=above right:$d_2$](d)at(1.2,4.4);

    \end{tikzpicture}
    \caption{$E_2$-term}
\end{figure}
\noindent
The bottom $d_2$ is zero similar to the case $d=2$. Note that $d_r(e^{*}_{i}\otimes e^{*}_{i})=0$ for arbitrary $r\geq 2$. This implies that each $e^{*}_{i}e^{*}_{j}\otimes e^{*}_{i}e^{*}_{j}$ is a permanent cycle, hence so is $(e^{*}_{i}e^{*}_{j}\otimes e^{*}_{i}e^{*}_{j})\alpha^2$. So $x\alpha^2$ cannot be hit. To vanish these three permanent cycles $(e^{*}_{i}e^{*}_{j}\otimes e^{*}_{i}e^{*}_{j})\alpha^2$, the rank of the $d_2$ at the top in the figure has to be $3$. By \cref{SWh}, the element $x$, given by the formula in the proof of \cref{confi_rep}, has to hit to $\alpha^4$ at $E_4$-term. So $d_2(x)$ and $d_3(x)$ are zero, and $d_4(x)=\alpha^4$. In summary, the rank of the top $d_2$ in the figure is $3$, and the others are zero. Therefore, the $E_3$-term is indicated as the following figure.
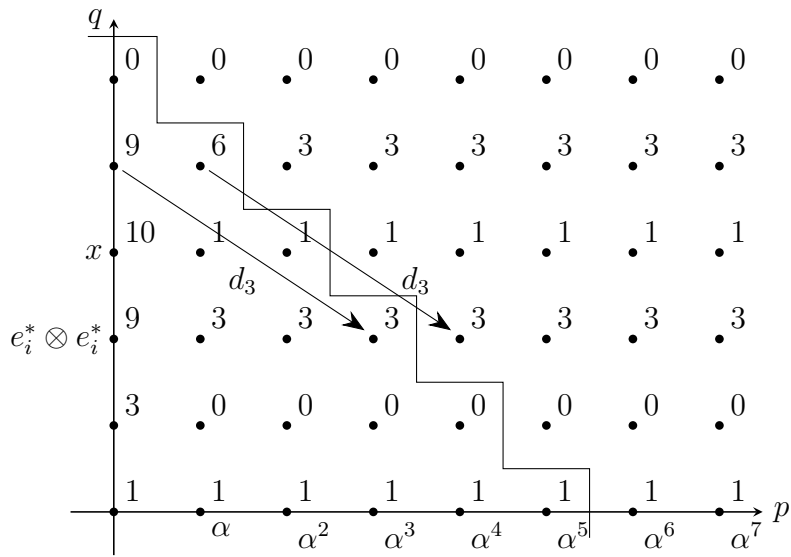
\begin{figure}[H]
    \centering
    \begin{tikzpicture}[xscale=1.15,yscale=1.15]
        \coordinate(O)at(0,0);
        \coordinate(XS)at(-0.5,0);
        \coordinate(XL)at(7.5,0);
        \coordinate(YS)at(0,-0.5);
        \coordinate(YL)at(0,5.7);
        \draw[semithick,->,>=stealth](XS)--(XL)node[right]{$p$};
        \draw[semithick,->,>=stealth](YS)--(YL)node[left]{$q$};
        \coordinate(P)at(1,0);
        \coordinate(Q)at(0,1);

        \foreach\k in{0,1,2,3,4,5,6,7}\foreach\l in{0,1,2,3,4,5}\fill($(O)+\k*(P)+\l*(Q)$)circle(0.05);
        
        \foreach\k in{0,1,2,3,4,5,6,7}
        \coordinate[label=above right:$1$](X\k)at($\k*(P)$);
        
        \coordinate[label=below right:$\alpha$](X)at(1,0);
        \foreach\k in{2,3,4,5,6,7}
        \coordinate[label=below right:$\alpha^\k$](X\k)at($\k*(P)$);
        
        \coordinate[label=above right:$3$](Y)at(0,1);
        
        \foreach\k in{1,2,3,4,5,6,7}
        \coordinate[label=above right:$0$](Y\k)at($\k*(P)+(Q)$);
        
        \coordinate[label=above right:$9$, label=left:$e^{*}_{i}\otimes e^{*}_{i}$](Z)at(0,2);
        \foreach\k in{1,2,3,4,5,6,7}\coordinate[label=above right:$3$](Z\k)at($\k*(P)+(Z)$);
        
        \coordinate[label=above right:$10$, label=left:$x$](W)at(0,3);
        
        \foreach\k in{1,2,3,4,5,6,7}\coordinate[label=above right:$1$](W\k)at($\k*(P)+(W)$);

        \coordinate[label=above right:$9$](R)at(0,4);
        
        \foreach\k in{1}\coordinate[label=above right:$6$](R\k)at($\k*(P)+(0,4)$);
        \foreach\k in{2,3,4,5,6,7}\coordinate[label=above right:$3$](R\k)at($\k*(P)+(0,4)$);

        \foreach\k in{0,1,2,3,4,5,6,7}\coordinate[label=above right:$0$](R\k)at($\k*(P)+(0,5)$);

        \draw(-0.3,5.5)--(0.5,5.5)--(0.5,4.5)--(1.5,4.5)--(1.5,3.5)--(2.5,3.5)--(2.5,2.5)--(3.5,2.5)--(3.5,1.5)--(4.5,1.5)--(4.5,0.5)--(5.5,0.5)--(5.5,-0.3);

        \draw [-{Stealth[length=3mm]}] (0.1,3.95) -- (2.9,2.1);
        \coordinate[label=above right:$d_3$](d)at(1.2,2.4);
        \draw [-{Stealth[length=3mm]}] (1.1,3.95) -- (3.9,2.1);
        \coordinate[label=above right:$d_3$](d)at(3.2,2.4);

    \end{tikzpicture}
    \caption{$E_3$-term}
\end{figure}
At the $E_3$-term, the remarkable differential is the right $d_3$ in the figure. Since its target is generated by the permanent cycles $(e^{*}_{i}\otimes e^{*}_{i})\alpha^4$, the rank of the right $d_3$ has to be $3$ to vanish such cycles at $E_\infty$. Hence, the left $d_3$ is rank $3$. Finally, $d_4(x)=\alpha^4$ implies that the $E_\infty$-term is indicated as the following figure.
\begin{figure}[H]
    \centering
    \begin{tikzpicture}[xscale=1.15,yscale=1.15]
        \coordinate(O)at(0,0);
        \coordinate(XS)at(-0.5,0);
        \coordinate(XL)at(7.5,0);
        \coordinate(YS)at(0,-0.5);
        \coordinate(YL)at(0,5.7);
        \draw[semithick,->,>=stealth](XS)--(XL)node[right]{$p$};
        \draw[semithick,->,>=stealth](YS)--(YL)node[left]{$q$};
        \coordinate(P)at(1,0);
        \coordinate(Q)at(0,1);

        \foreach\k in{0,1,2,3,4,5,6,7}\foreach\l in{0,1,2,3,4,5}\fill($(O)+\k*(P)+\l*(Q)$)circle(0.05);
        
        \foreach\k in{0,1,2,3}
        \coordinate[label=above right:$1$](X\k)at($\k*(P)$);
        \foreach\k in{4,5,6,7}
        \coordinate[label=above right:$0$](X\k)at($\k*(P)$);
        
        \coordinate[label=below right:$\alpha$](X)at(1,0);
        \foreach\k in{2,3}
        \coordinate[label=below right:$\alpha^\k$](X\k)at($\k*(P)$);

        \coordinate[label=above right:$3$](Y)at(0,1);
        
        \foreach\k in{1,2,3,4,5,6,7}
        \coordinate[label=above right:$0$](Y\k)at($\k*(P)+(Q)$);
        
        \coordinate[label=above right:$9$](Z)at(0,2);
        \foreach\k in{1,2}\coordinate[label=above right:$3$](Z\k)at($\k*(P)+(Z)$);
        \foreach\k in{3,4,5,6,7}\coordinate[label=above right:$0$](Z\k)at($\k*(P)+(Z)$);
        
        \coordinate[label=above right:$9$](W)at(0,3);
        
        \foreach\k in{1,2,3,4,5,6,7}\coordinate[label=above right:$0$](W\k)at($\k*(P)+(W)$);

        \coordinate[label=above right:$6$](R)at(0,4);
        
        \foreach\k in{1}\coordinate[label=above right:$3$](R\k)at($\k*(P)+(0,4)$);
        \foreach\k in{2,3,4,5,6,7}\coordinate[label=above right:$0$](R\k)at($\k*(P)+(0,4)$);

        \foreach\k in{0,1,2,3,4,5,6,7}\coordinate[label=above right:$0$](R\k)at($\k*(P)+(0,5)$);

        \draw(-0.3,5.5)--(0.5,5.5)--(0.5,4.5)--(1.5,4.5)--(1.5,3.5)--(2.5,3.5)--(2.5,2.5)--(3.5,2.5)--(3.5,1.5)--(4.5,1.5)--(4.5,0.5)--(5.5,0.5)--(5.5,-0.3);

    \end{tikzpicture}
    \caption{$E_\infty$-term}
\end{figure}
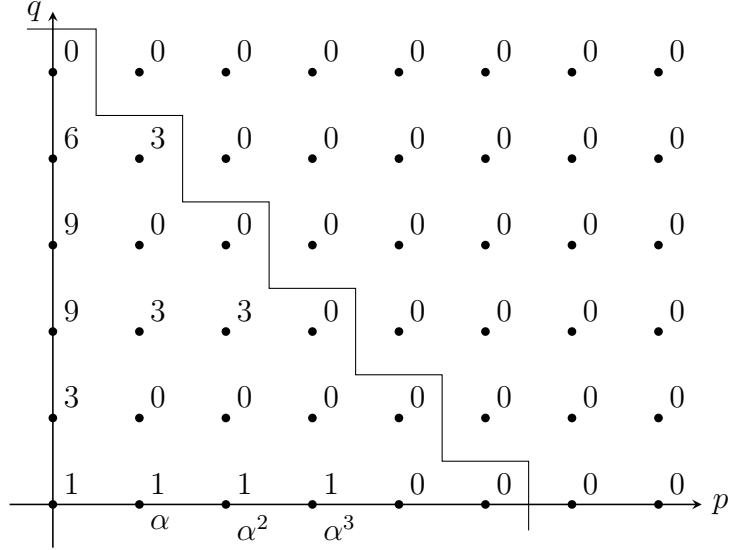

Finally, we choose generators $u_{i,j}$ and $x_{i,j}$ as generators of $E_{\infty}^{0,i}$. In particular, let $u_{i,j}$ be in the regular representation parts, and $x_{i,j}$ in the trivial representation parts. Here we pay attention to that each $u_{i,j}\alpha^k(k\geq 1)$ might be nonzero in $\ccohomology{*}{\uconfi{2}{\torus{d}}}{\F_2}$ in spite of it is zero in $E_{\infty}$ (extension problem). However, we can choose $u_{i,j}$ so that $u_{i,j}\alpha$ is zero in $\ccohomology{*}{\uconfi{2}{\torus{d}}}{\F_2}$. Indeed, if $u_{i,j}\alpha$ is nonzero in $\ccohomology{*}{\uconfi{2}{\torus{d}}}{\F_2}$, there is a positive integer $m$ and an element $y\in \ccohomology{i}{\uconfi{2}{\torus{d}}}{\F_2}$ such that $u_{i,j}\alpha=y\alpha^m$. Then replacing $u_{i,j}$ to $u'_{i,j}=u_{i,j}+y\alpha^{m-1}$, we obtain $u'_{i,j}\alpha=0$ in $\ccohomology{*}{\uconfi{2}{\torus{d}}}{\F_2}$. This completes the proof of \cref{case_2_3}.

\section*{Acknowledgement}
The author appreciates Mitsunobu Tsutaya and Takahiro Matsushita for commenting on this study. The author was also supported by WISE program (MEXT) at Kyushu University.

\vspace{3mm}

\end{document}